\documentclass[12pt,oneside]{article}
\usepackage{amsthm,amsmath,amssymb,txfonts,graphics,,mathrsfs}
\usepackage[latin1]{inputenc}

\usepackage{pb-diagram}

\newtheorem{defin}{Definition}
\newtheorem{lemma}{Lemma}
\newtheorem{teo}{Theorem}
\newtheorem{oss}{Remark}
\newtheorem{prop}{Proposition}
\newtheorem{coroll}{Corollary}

\author{Jung Kyu {CANCI}}
\title{Cycles for rational maps of good reduction outside a prescribed set
}
\date{}

\begin{document}
\maketitle
\begin{abstract}\noindent Let $K$ be a number field and $S$ a fixed finite set of places of $K$ containing all the archimedean ones. Let $R_S$ be the ring of $S$-integers of $K$. In the present paper we study the cycles in $\mathbb{P}_1(K)$ for rational maps of degree $\geq2$ with good reduction outside $S$. We say that two ordered $n$-tuples $(P_0,P_1,\ldots,P_{n-1})$ and $(Q_0,Q_1,\ldots,Q_{n-1})$ of points of $\mathbb{P}_1(K)$ are equivalent if there exists an automorphism $A\in{\rm PGL}_2(R_S)$ such that $P_i=A(Q_i)$ for every index $i\in\{0,1,\ldots,n-1\}$. We prove that if we fix two points $P_0,P_1\in\mathbb{P}_1(K)$, then the number of inequivalent cycles for rational maps of degree $\geq2$ with good reduction outside  $S$ which admit $P_0,P_1$ as consecutive points is finite and depends only on $S$ and $K$. We also prove that this result is in a sense best possible.
\end{abstract}
\section{Introduction}
\noindent Let $K$ be a number field and $R$ its ring of integers. Let $\Phi\colon\mathbb{P}_1\to\mathbb{P}_1$ be a rational map defined over $K$ by $\Phi([x:y])=\left[F(x,y):G(x,y)\right]$ where $F,G\in R\left[x,y\right]$ are homogeneous of the same degree with no common factor.

Let $\mathfrak{p}$ be a prime ideal of $R$. Using the standard notation which will be introduced at the beginning of next section, we can assume that $F,G$ have coefficients in $R_{\mathfrak{p}}$ (the local ring of $R$ at the prime ideal $\mathfrak{p}$) and at least one coefficient belonging to $R_{\mathfrak{p}}^*$. In this way, let $K(\mathfrak{p})=R/\mathfrak{p}$ be the residue field, we obtain the rational map defined over $K(\mathfrak{p})$
\begin{displaymath} \tilde{\Phi}\colon\mathbb{P}_1\to\mathbb{P}_1;\ \ \tilde{\Phi}([x:y])=[\tilde{F}(x,y):\tilde{G}(x,y)]\end{displaymath} 
where $\tilde{F},\tilde{G}$ are the polynomials obtained by reducing modulo $\mathfrak{p}$ the coefficients of $F$ and $G$. Since the polynomial $F,G$ have no common factors, the rational map $\tilde{\Phi}$ has degree equal to $\deg \Phi$ if and only if $\tilde{F}$ and $\tilde{G}$ do not have common roots in $\mathbb{P}_1(\overline{K(\mathfrak{p})})$, where $\overline{K(\mathfrak{p})}$ denotes the algebraic closure of $K(\mathfrak{p})$. If this is the case we will say that $\Phi$ has good reduction at the prime ideal $\mathfrak{p}$.

A cycle of length $n$ for a rational map $\Phi$ is a ordered $n$--tuple $(P_0, P_1,\ldots,P_{n-1})$, of distinct points of  $\mathbb{P}_1(K)$, with the property that $\Phi(P_i)=P_{i+1}$ for every 
$i\in\{0,1,
\ldots,n-2\}$ and such that $\Phi(P_{n-1})=P_0$. It is easy to see that every $n$--tuple of distinct points is a cycle for a suitable rational map, but imposing some restrictions on the maps gives some considerable restrictions on the cycles. In the present paper we shall consider a number field $K$ and a fixed finite set $S$ of places containing the archimedean ones. We will study the cycles for rational maps that have good reduction at every prime ideal $\mathfrak{p}$ whose associated $\mathfrak{p}$-adic place does not belong to  $S$ and we will say that these maps \emph{have good reduction outside} $S$.

Morton and Silverman \cite[Corollary B]{M.S.1} have proved that if $\Phi$ is a rational map of degree $\geq 2$ which has bad reduction only at $s$ prime ideals of $K$ and $P\in\mathbb{P}_1(K)$ is a periodic point with minimal period $n$, then the following inequality holds:
\begin{displaymath} n\leq (12(s+2)\log(5(s+2)))^{4\left[K:\mathbb{Q}\right]}.\end{displaymath} 
These results provide some bounds for the period-length of a periodic point for a rational map $\Phi$ depending only on the number of prime ideals of bad reduction.
This generalizes and improves the result by Narkiewicz \cite{N.1} who was concerned with polynomial maps, in fact if $\phi(z)\in K[z]$ is a polynomial, then the corresponding map $\phi\colon\mathbb{P}_1\to\mathbb{P}_1$ has good reduction outside $S$ if and only if $\phi$ has $S$-integral coefficients and its leading coefficient is a $S$-unit. For this type of polynomials, Narkiewicz found a bound for minimal period-length $n$ which is possible to write in the following way:
$n\leq C^{s^2+s\left[K:\mathbb{Q}\right]}$, where $C$ is an absolute constant.\\
The main tool used by Narkiewicz is the finiteness of the solutions in $S$-units of the equation $u+v=1$, and in particular the estimate of Evertse \cite{E.1}. On the other hand, Morton and Silverman used their results on multiplicity and reduction obtained in \cite{M.S.2}.
R. Benedetto has recently obtained a much stronger bound for polynomial maps. He proved in \cite{B.1} that if $\phi\in K[z]$ is a polynomial of degree $d\geq2$ which has bad reduction in $s$ primes of $K$, then the number of preperiodic points of $\phi$ is at most $O(s\log s)$. The big-$O$ constant is essentially $(d^2-2d+2)/\log d$ for large $s$. Benedetto's proof relies on a detailed analysis of $\mathfrak{p}$-adic Julia sets.

Let $R_S$ be the ring of $S$-integers of $K$; the automorphism-group ${\rm PGL}_2(R_S)$ acts in a canonical way on $\mathbb{P}_1(K)$. If $(P_0, P_1,\ldots,P_{n-1})$ is a $n$-cycle for a rational map $\Phi$, with good reduction outside $S$, then for every $A\in PGL_2(R_S)$ the \emph{image}-$n$-tuple  $(A(P_0),A(P_1),\ldots,A(P_{n-1}))$ is a  $n$-cycle for the rational map $A\circ\Phi \circ A^{-1}$, which still has good reduction outside $S$; we will call the two $n$-tuples \emph{equivalent}.

In \cite[Theorem 1]{H.N.1} Halter-Koch and Narkiewicz proved the finiteness of the set of possible normalized $n$-cycles in $R_S$ for polynomial maps, where a cycle is called normalized when $0$ and $1$ are two consecutive elements of the $n$-tuple. In the present paper we generalize to rational maps this result, in particular we prove the following corollary as a consequence of our Theorem 1 below.

\begin{coroll}\label{P_0,P_1}
Let $P_0,P_1\in\mathbb{P}_1(K)$ be two fixed points. The number of inequivalent cycles for rational maps of degree $\geq2$ with good reduction outside  $S$ which admit $P_0,P_1$ as consecutive points is finite and depends only on $S$ and $K$.\end{coroll}

Let $P_1=\left[x_1:y_1\right],P_2=\left[x_2:y_2\right]\in \mathbb{P}_1(K)$ and $\mathfrak{p}$ a prime ideal of $R$. Using the notation of  \cite{M.S.2} we will denote by  \begin{equation}\label{d_p}\delta_{\mathfrak{p}}\,(P_1,P_2)=v_{\mathfrak{p}}\,(x_1y_2-x_2y_1)-\min\{v_{\mathfrak{p}}(x_1),v_{\mathfrak{p}}(y_1)\}-\min\{v_{\mathfrak{p}}(x_2),v_{\mathfrak{p}}(y_2)\}\end{equation}the $\mathfrak{p}$-adic logarithmic distance; $\delta_{\mathfrak{p}}\,(P_1,P_2)$ is independent of the choice of the homogeneous coordinates, i.e. it is well defined.
\\To every pair $P,Q\in\mathbb{P}_1(K)$ we associate the ideal
\begin{displaymath} \mathfrak{I}(P,Q):=\prod_{\substack{\mathfrak{p}\notin S}}\mathfrak{p}^{\delta_{\mathfrak{p}}(P,Q)}.\end{displaymath}

It is characterized by the property that $P\equiv Q\ ({\rm mod}\ \mathfrak{I}(P,Q))$ and that for every ideal $\mathfrak{I}$ such that $P\equiv Q\ ({\rm mod}\ \mathfrak{I})$ one has $\mathfrak{I}\mid\mathfrak{I}(P,Q)$.\\
To every $n$-tuple $(P_0, P_1,\ldots,P_{n-1})$ we can associate the $(n-1)$-tuple of ideals $(\mathfrak{I}_1,\mathfrak{I}_2,\ldots,\mathfrak{I}_{n-1})$ defined by  
\begin{equation}\label{I_i} \mathfrak{I}_i:=\prod_{\substack{\mathfrak{p}\notin S}}\mathfrak{p}^{\delta_{\mathfrak{p}}(P_0,P_i)}=\mathfrak{I}(P_0,P_i)\end{equation}
\\
With the above notation we will prove the following results.
 
\begin{teo}\label{tp} There exists a finite set $\mathbb{I}_S$ of ideals of $R_S$, depending only on $S$ and $K$, with the following property: for every $n$-cycle $(P_0, P_1,\ldots,P_{n-1})$, for a rational map of degree $\geq2$ with good reduction outside $S$, let $(\mathfrak{I}_1,\mathfrak{I}_2,\ldots,\mathfrak{I}_{n-1})$ be the associated $(n-1)$-tuple of ideals; then \begin{displaymath} \mathfrak{I}_i\mathfrak{I}_1^{-1}\in\mathbb{I}_S\end{displaymath} 
for every index $i\in\{1,\ldots,n-1\}$.\end{teo}
Theorem \ref{tp} will be first proved with the particular condition that $R_S$ is a P.I.D.. Afterwards the proof in the general case will follow.

The proof of Corollary \ref{P_0,P_1} will be a direct consequence of Theorem \ref{tp} by applying the results obtained by Birch and Merriman in 1972 \cite{B.M.1}.

In the proof of Theorem \ref{tp} we will also prove 
\begin{coroll}\label{N}
There exists a finite set $\mathcal{N}$ of $n$-tuples depending only on $S$ and $K$, such that if $\Phi$ is a rational map of degree $\geq 2$ with good reduction outside $S$, then every $n$-cycle for $\Phi$ in $\mathbb{P}_1(K)$ can be transformed by an automorphism in ${\rm PGL}_2(K)$ into an $n$-tuple in $\mathcal{N}$. 
\end{coroll}

Two cycles $(x_0,x_1,\ldots,x_{n-1})$ and $(y_0,y_1,\ldots,y_{n-1})$ for polynomial maps are called equivalent if and only if there exist an $S$-integer $a\in R_S$ and an $S$-unit $\epsilon\in R_S^\ast$ such that $y_i=a+\epsilon x_i$, for every index $i$. The definition of equivalent cycles for rational maps introduced above is the natural generalization of the one just stated.

Theorem 2 in \cite{H.N.1} states that in $R_S$ for every $n\geq2$ there is just a finite number of inequivalent $n$-cycles for polynomial maps of degree $\geq2$.
This result cannot be extended to rational maps, since we have proved the following theorem in the case $S$ contains the places which extend the $2$-adic place of $\mathbb{Q}$.

\begin{teo}\label{tp2} Let $R_S=\mathbb{Z}[1/2]$. There exist infinitely many ideals $\mathfrak{I}$ for which there exists a $3$-cycle $(P_0,P_1,P_2)$, for a suitable rational map of degree $4$ with good reduction outside $S$, for which $\mathfrak{I}_1=\mathfrak{I}$ holds, where $\mathfrak{I}_1$ is the ideal defined in (\ref{I_i}).
\end{teo}

Theorem \ref{tp2} proves that the conclusion of Theorem \ref{tp} is in a sense best-possible: for every cycle one has
\begin{displaymath} (\mathfrak{I}_1,\ldots,\mathfrak{I}_{n-1})=\mathfrak{I}_1(R_S,\mathfrak{I}_2\mathfrak{I}_1^{-1},\ldots,\mathfrak{I}_{n-1}\mathfrak{I}_1^{-1})\end{displaymath}
where for the factor $(R_S,\mathfrak{I}_2\mathfrak{I}_1^{-1},\ldots,\mathfrak{I}_{n-1}\mathfrak{I}_1^{-1})$ there are only finitely many possibilities  in view of Theorem \ref{tp}, but not for the factor $\mathfrak{I}_1$, in view of Theorem \ref{tp2}.

Our method of proof is similar to the one used by W. Narkiewicz, F. Halter-Koch and T. Pezda (see \cite{H.N.1}, \cite{N.1}, \cite{N.P.1}, \cite{P.1}). It provides an effective bound for the cardinality of $\mathbb{I}_S$ depending only on $S$. Unfortunately, by the same method, we cannot obtain an effective solution to the problem, since we shall use the theorem on the  finiteness of solutions to equations in three $S$-units $u_1,u_2,u_3$: 
\begin{displaymath} a_1u_1+a_2u_2+a_3u_3=1,\end{displaymath}
which is not effective. 

\medskip
\noindent$\mathbf{Acknowledgments.}$ The present work will be part of my Ph.D.-thesis supervised by professor P. Corvaja. I would like to thank him for his useful advice.
I would like to thank also professor Narkiewicz who has sent me his fundamental articles concerning
this topic. I am grateful to the referee for correcting some inaccuracies on a previous draft of the present paper and for suggesting substantial simplifications.

\section{Good reduction for $n$-tuples}
In all the present paper we will use the following notation:
\begin{itemize}
\item[$K$]a number field;
\item[$R$]the ring of integers of $K$;
\item[$\mathfrak{p}$]a prime ideal of $R$, $\mathfrak{p}\neq0$;
\item[$R_\mathfrak{p}$]the local ring of $R$ at the prime ideal $\mathfrak{p}$;
\item[$m_\mathfrak{p}$] the maximal ideal of $R_\mathfrak{p}$ (which is principal);
\item[$K(\mathfrak{p})$]=$R/\mathfrak{p}\cong R_{\mathfrak{p}}/m_{\mathfrak{p}}$ the residue field of the prime ideal $\mathfrak{p}$;
\item[$v_{\mathfrak{p}}$]the $\mathfrak{p}$-adic valuation on $R$ corresponding to the prime ideal $\mathfrak{p}$ (we always assume $v_\mathfrak{p}$ to be normalized so that $v_{\mathfrak{p}}(K^*)=\mathbb{Z}$); 
\item[$S$]a fixed finite set of places of $K$ of cardinality $s$ including all archimedean places.
\end{itemize}
We denote the ring of $S$-integers by
\begin{displaymath}  R_S:=\{x\in K \mid v_{\mathfrak{p}}(x)\geq0 \ \text{for every prime ideal }\ \mathfrak{p}\notin S\}
\end{displaymath}
and the group of $S$-units by
\begin{displaymath}  R_S^\ast:=\{x\in K^\ast\mid v_{\mathfrak{p}}(x)=0 \ \text{for every prime ideal }\ \mathfrak{p}\notin S\}.
\end{displaymath}

The canonical $({\rm mod }\ \mathfrak{p})$-projection $\mathbb{P}_1(K)\to \mathbb{P}_1(K(\mathfrak{p}))$ is defined in the following way: for every point $P\in \mathbb{P}_1(K)$, choose some coordinates $P=\left[x:y\right]$ such that $x,y\in R_{\mathfrak{p}}$ and they do not belong simultaneously to $m_{\mathfrak{p}}$, so the point $\left[x+m_{\mathfrak{p}}:y+m_{\mathfrak{p}}\right]\in \mathbb{P}_1(R_{\mathfrak{p}}/m_{\mathfrak{p}})$ is well defined. By the canonical isomorphism $R_{\mathfrak{p}}/m_{\mathfrak{p}}\cong K(\mathfrak{p})$, for every point $P\in \mathbb{P}_1(K)$ it is possible to associate a point of $\mathbb{P}_1(K(\mathfrak{p}))$ which will be called the \emph{reduction modulo $\mathfrak{p}$} of $P$. 
\begin{defin}\label{nbr}
Let $\mathfrak{p}\neq0$ be a prime ideal of $R$. We say that a $n$-tuple $(P_0,\ldots,P_{n-1})$ of elements of $\ \mathbb{P}_1(K)$ has good reduction at $\mathfrak{p}$ if the $n$-tuple formed by the reduction modulo $\mathfrak{p}$ has $n$ distinct elements of $\ \mathbb{P}_1(K(\mathfrak{p}))$; a $n$-tuple has good reduction outside $S$ if it has good reduction at every prime ideal $\mathfrak{p}\notin S$.
\end{defin}
The finiteness of the class number of the ring $R$ will be used to prove

\begin{prop}\label{cqc}
There exists a finite set $S_R$ of non-archimedean places of $K$ and an integer $C$, depending only on $R$, such that every point $P\in\mathbb{P}_1(K)$ can be represented by integral homogeneous coordinates $(x,y)$ satisfying $\min\{v_{\mathfrak{p}}(x),v_{\mathfrak{p}}(y)\}=0$ for all prime ideal $\mathfrak{p}\notin S_R$ and $\min\{v_{\mathfrak{p}}(x),v_{\mathfrak{p}}(y)\}\leq C$ for every $\mathfrak{p}\in S_R$ .\end{prop}
\begin{proof} Let \mbox{$R,(a_2R+b_2R),\ldots,(a_tR+b_tR)$} be a set of representatives for the ideal classes.
Each point $P\in\mathbb{P}_1(K)$ can be expressed by integer coordinates $P=[x:y]$ such that $(xR+yR)=(a_iR+b_iR)$ for a suitable index $i\in\{1,2,\ldots,t\}$ (with $a_1=b_1=1$), thus 
\begin{displaymath} \min\{v_{\mathfrak{p}}(x),v_{\mathfrak{p}}(y)\}=\min\{v_{\mathfrak{p}}(a_i),v_{\mathfrak{p}}(b_i)\}\end{displaymath}
for every prime ideal $\mathfrak{p}$.\\
Now it is sufficient to choose $S_R$ as the set of non-archimedean places such that $\min\{v_{\mathfrak{p}}(a_i),v_{\mathfrak{p}}(b_i)\}\neq 0$ for some index $i\in\{2,\ldots,t\}$ and \begin{equation}\label{log} C=\max_{\mathfrak{p}\in S_R}\left\{\min\{v_{\mathfrak{p}}(a_i),v_{\mathfrak{p}}(b_i)\}\mid i\in\{2,\ldots,t\}\right\}.\end{equation} 
The constant $C$ and the set $S_R$ can be taken depending only on $K$ since \cite[Corollary 2 to Theorem 36, Chapter 5]{M.1}
\end{proof}

Proposition \ref{cqc} allows to adopt the following convention: writing \mbox{$P=\left[x:y\right]$} for a generic element of $\mathbb{P}_1(K)$ we will always choose $x,y\in R$  with the property just described and we will say that $x$ and $y$ are \emph{almost coprime}.

\vspace{3.5mm}\noindent$\mathbf{Notation.}$ In the present section every point will be represented with almost coprime coordinates, except in the cases in which it will be explicitly specified. Moreover for any $n$-tuple $(P_0, P_1,\ldots,P_{n-1})$ of points of $\mathbb{P}_1(K)$, for every index $i$, $(x_i,y_i)$ always will represent almost coprime integral homogeneous coordinate for the point $P_i$.

\vspace{3.5mm}The $\mathfrak{p}$-adic logarithmic distance $\delta_{\mathfrak{p}}$ defined in (\ref{d_p}) assumes integral values and the following properties hold:
\begin{displaymath} \left.\begin{array}{lll}
\delta_{\mathfrak{p}}(P,Q)\geq0\ & \text{for every $P$ and $Q$} & \textrm{($\delta^\prime$)}\\
\delta_{\mathfrak{p}}(P,Q)\geq1\ & \text{if and only if}\  P\equiv Q\ (\rm{mod}\ \mathfrak{p}) & \textrm{($\delta^{\prime\prime}$)}\\
\delta_{\mathfrak{p}}(P,Q)=\infty\ & \text{if and only if}\  P=Q & \textrm{($\delta^{\prime\prime\prime}$)}
\end{array}\right.\end{displaymath}
By property ($\delta^{\prime\prime}$) it follows that a $n$-tuple $(P_0,P_1,\ldots,P_{n-1})\in\mathbb{P}_1^n(K)$ has good reduction outside $S$ if and only if $\delta_{\mathfrak{p}}\,(P_i,P_j)=0$ for every  prime ideal $\mathfrak{p}$ not in $S$ and for every distinct indexes $i,j\in\{0,\ldots,n-1\}$. \\
Therefore if the $n$-tuple $(P_0, P_1,\ldots,P_{n-1})$ has good reduction outside $S$, then the $(n-1)$-tuple $(\mathfrak{I}_1,\mathfrak{I}_2,\ldots,\mathfrak{I}_{n-1})$ of ideals defined by (\ref{I_i}) is equal to $(R_S,\ldots,R_S)$.
\begin{defin}\label{d3}
Two $n$-tuples $(P_0,P_1,\dots,P_{n-1})$ and $(Q_0,Q_1,\dots,Q_{n-1})$ are called equivalent if there exists a projective automorphism $A \in {\rm PGL}_2(R_S)$ such that 
\begin{displaymath}A(P_i)=Q_i\  for\  all\  i\in\{0,1,\dots,n-1\}.
\end{displaymath}
\end{defin}
If the $n$-tuples $(P_0,\dots,P_{n-1})$ and $(Q_0,\dots,Q_{n-1})$ are equivalent, then the $(n-1)$-tuples $(\mathfrak{I}_1,\mathfrak{I}_2,\ldots,\mathfrak{I}_{n-1})$ of ideals defined by (\ref{I_i}) coincide. Moreover if a $n$-tuple $(P_0, P_1,\ldots,P_{n-1})$ has good reduction outside $S$, then every $n$-tuple equivalent to $(P_0, P_1,\ldots,P_{n-1})$  has good reduction outside $S$ as well.

If $R_S$ is a P.I.D. (principal ideal domain), then the class number of $R_S$ is $1$ and $R_S$ is a representative of the unique ideal class. In this case, for any point of $\mathbb{P}_1(K)$ we can choose coprime integral homogeneous coordinates and then, from this choice, for any $n$-tuple $(P_0, P_1,\ldots,P_{n-1})$ with good reduction outside $S$ and for every prime ideal $\mathfrak{p}\notin S$ it follows that 
\begin{displaymath} v_{\mathfrak{p}}(x_iy_j-x_jy_i)=\delta_{\mathfrak{p}}(P_i,P_j)=0\end{displaymath}
thus $x_iy_j-x_jy_i$ is a $S$-unit.\\
In general $R_S$ is not always a P.I.D.. Enlarging $S$ (to obtain a principal domain) we change the equivalence relation between $n$-tuples. So we can not change $S$ since we investigate about finiteness of the orbits of $n$-tuples under the action of ${\rm PGL}_2(R_S)$. But in any case we have that

\begin{lemma}\label{obrn}
There exists a finite set $\mathcal{R}$ of $S$-integers depending only on $S$ and $K$ such that for any $n$-tuple $(P_0, P_1,\ldots,P_{n-1})$  of good reduction outside $S$ and for every distinct indexes $i,j\in\{0,\ldots,n-1\}$ there exist $r_{i,j}\in \mathcal{R}$ and a unit $u_{i,j}\in R_S^\ast$ such that
\begin{displaymath} x_iy_j-x_jy_i=r_{i,j}u_{i,j}.\end{displaymath} \end{lemma}
\begin{proof} By the good reduction of $(P_0, P_1,\ldots,P_{n-1})$ and the definition of logarithmic distance we have that 
\begin{displaymath} v_{\mathfrak{p}}(x_iy_j-x_jy_i)=\min\{v_{\mathfrak{p}}(x_i),v_{\mathfrak{p}}(y_i)\}+\min\{v_{\mathfrak{p}}(x_j),v_{\mathfrak{p}}(y_j)\}\end{displaymath}
for every $\mathfrak{p}\notin S$.
Let $C$ be the integer and $S_R$ the set defined in the Proposition \ref{cqc}. Having chosen almost coprime coordinates it follows that 
\begin{equation}\label{2C} v_{\mathfrak{p}}(x_iy_j-x_jy_i)\leq 2C,\end{equation} for every $\mathfrak{p}\in S_R/S$, and $v_{\mathfrak{p}}(x_iy_j-x_jy_i)=0$ for every other prime ideal not in $S$.\\
For every couple of distinct points $P_i=[x_i:y_i],P_j=[x_j:y_j]$ included in a $n$-tuple which has good reduction outside $S$, it is defined the following principal ideal of $R_S$ 
\begin{displaymath} (x_iy_j-x_jy_i)R_S=\prod_{\substack{\mathfrak{p}\in S_R\setminus S}}\mathfrak{p}^{v_{\mathfrak{p}}(x_iy_j-x_jy_i)}. \end{displaymath} 
By (\ref{2C}) we are in the position to conclude that the set of principal ideals generated in this way has finite cardinality and therefore choosing a generator for every ideal defines the finite set $\mathcal{R}$ which has cardinality bounded by $(2C+1)^{|S_R\setminus S|}$.
By the remarks made at the end of the proof of Proposition \ref{cqc} we deduce that the cardinality of $\mathcal{R}$ is bounded by a constant which depends only on $S$ and $K$.\end{proof}

By the finiteness of classes of binary forms with given discriminant proved by Birch and Merriman in 1972 \cite{B.M.1} we deduce the following
\begin{prop}\label{fn}
The set of equivalence classes of $n$-tuples with good reduction outside $S$ is finite and depends only on the set $S$ and $K$. 
\end{prop}
\begin{proof}
Note first that for large $n$ there are no $n$-tuples with good reduction outside $S$. Indeed,
let $m=\min_{\mathfrak{p}\notin S}\{\left|K(\mathfrak{p})\right|\}$. For every $n\geq m+2$, each $n$-tuple will not be of good reduction outside $S$, since for every prime ideal $\mathfrak{p}$ which realizes the minimum $m$ the projective space $\mathbb{P}_1(K(\mathfrak{p}))$ has only $m+1$ elements. \\
For $n=1$ the number of equivalence classes is the order of the ideal class group. Indeed, let $r$ be the class number of $R_S$. We choose the representatives for every class and express them, except the trivial ideal $R_S$, through two generators $(a_2R_S+b_2R_S),\ldots,(a_rR_S+b_rR_S)$. Note that these representatives define $r$ points of $\mathbb{P}_1(K)$ inequivalent for the action of ${\rm PGL}_2(R_S)$. Now we prove that every point $P\in\mathbb{P}_1(K)$ belongs to the orbit of [1:0] or of a point $[a_i:b_i]$ for a suitable index $i\in\{2,\ldots,r\}$.\\
Let $P\in\mathbb{P}_1(K)$. We write it with integer coordinates $P=[\bar{x}:\bar{y}]$. \\
If $(\bar{x}R_S+\bar{y}R_S)$ is a principal ideal (i.e. it is equivalent to the trivial ideal $R_S$), then $P$ is an element of the orbit of $[1:0]$ under the action of ${\rm PGL}_2(R_S)$. Indeed, let $(\bar{x}R_S+\bar{y}R_S)=aR_S$ for a suitable $a\in R_S$; then $x=\bar{x}/a$ and $y=\bar{y}/a$ are elements of $R_S$ such that $(xR_S+yR_S)=R_S$; this is equivalent to the existence of two $S$-integers $r_x$ and $r_y$ such that $xr_x+yr_y=1$; therefore the matrix $\begin{pmatrix}
	r_x & r_y\\
	-y & x
\end{pmatrix}$ belongs to ${\rm SL}_2(R_S)$
and maps the point $[x:y]$ to $[1:0]$.\\
Otherwise, there exist $c,d\in R_S$ such that $c(\bar{x}R_S+\bar{y}R_S)=d(a_iR_S+b_iR_S)$ for some index $i$; therefore, denoting $x=c\bar{x}/d$ and $y=c\bar{y}/d$ (elements of $R_S$), the following holds
\begin{displaymath} (xR_S+yR_S)=(a_iR_S+b_iR_S)=I\subset R_S. \end{displaymath}
By definition of $I^{-1}$, there are elements $x^{\prime},y^{\prime}\in I^{-1}$ satisfying $xy^{\prime}-yx^{\prime}=1$, namely $\begin{pmatrix}
	x & x^{\prime}\\
	y & y^{\prime}
\end{pmatrix}\in {\rm SL}_2(K)$. Moreover there are $a_i^{\prime},b_i^{\prime}\in I^{-1}$ such that $a_ib_i^{\prime}-b_ia_i^{\prime}=1$, namely $\begin{pmatrix}
	a_i & a_i^{\prime}\\
	b_i & b_i^{\prime}
\end{pmatrix}\in {\rm SL}_2(K)$. So that the following matrix
\begin{multline*} \begin{pmatrix}
	a_i & a_i^{\prime}\\
	b_i & b_i^{\prime}\end{pmatrix}\begin{pmatrix}
	x & x^{\prime}\\
	y & y^{\prime}
\end{pmatrix}^{-1}=\begin{pmatrix}
	a_i & a_i^{\prime}\\
	b_i & b_i^{\prime}\end{pmatrix}\begin{pmatrix}
	y^{\prime} & -x^{\prime}\\
	-y & x
\end{pmatrix}=\begin{pmatrix}
	a_iy^{\prime}-ya_i^{\prime} & -a_ix^{\prime}+xa_i^{\prime}\\
	y^{\prime}b_i-yb_i^{\prime} & -x^{\prime}b_i+xb_i^{\prime}
\end{pmatrix}\in {\rm SL}_2(R_S)
\end{multline*}
maps $[x:y]$ to $[a_i:b_i]$. This concludes the case $n=1$.\\
Let $n\geq2$ and $(P_0,P_1,\dots,P_{n-1})$ be a $n$-tuple with good reduction outside $S$. By Lemma \ref{obrn} we obtain for every distinct indexes $i,j$ the following identity
\begin{equation}\label{rbr} x_iy_j-x_jy_i=r_{i,j}u_{i,j},\end{equation} 
where $u_{i,j}\in R_S^*$ and $r_{i,j}\in\mathcal{R}$. For every choice of a unit $\lambda\in R_S^\ast$, these identities are verified also after replacing the almost coprime coordinates $\left[x_i:y_i\right]$ with the coordinates $\left[\lambda x_i:\lambda y_i\right]$ .
Now fixing the possible values of $r_{i,j}$, to every $n$-tuple of points which verifies the identities (\ref{rbr}) we associate the following binary form of degree $n$
\begin{displaymath} F(X,Y)=\prod_{0\leq i\leq n-1}{(x_iX-y_iY)},\end{displaymath}
defining in this way a family of forms with discriminant
\begin{equation}\label{D(F)} D(F)=u\left(\prod_{\substack{0\leq i<j\leq n-1}}{r_{i,j}^2}\right),\end{equation} 
where $u$ is a $S$-unit.
 
The multiplicative group $R_S^*$ is finitely generated, so there exists a finite set $\mathcal{V}\subset R_S^\ast$ such that every $S$-unit is representable as product of a $(2n-2)$-power of an $S$-unit and an element of $\mathcal{V}$. Then, for every $n$-tuple $(P_0, P_1,\ldots,P_{n-1})$ with good reduction which satisfies the identities (\ref{rbr}), the unit $u$ which appears in the equation (\ref{D(F)}) relative to the discriminant of the binary form associated to $(P_0, P_1,\ldots,P_{n-1})$ can be written as  $u=v\lambda^{2n-2}$ with $v\in\mathcal{V}$ and $\lambda\in R_S^\ast$. Thus, if we replace the coordinates $(x_0,y_0)$ of $P_0$ with $(x_0\lambda^{-1},y_0\lambda^{-1})$, we obtain a new binary form with discriminant $v\prod{r_{i,j}^2}$. \\
In other words, with an appropriate choice of coordinates, the $n$-tuples which satisfy the identities (\ref{rbr}) with fixed $r_{i,j}$ define a family of binary forms of degree $n$ whose discriminant is of the form $v\prod{r_{i,j}^2}$ with $v\in\mathcal{V}$.

Two binary forms $G(X,Y)$ and $H(X,Y)$ are called equivalent if there exists a matrix $A=\begin{pmatrix}
	a & b\\
	c & d
\end{pmatrix}\in {\rm GL}_2(R_S)$ such that $G(X,Y)=H(aX+bY,cX+dY)$. The equivalence of binary forms associated with $n$--tuples coincide with the equivalence of the corresponding unordered $n$--tuples.
By the results obtained by Birch and Merriman in 1972 \cite{B.M.1} (non effective) and Evertse and Gy\H ory in 1991 \cite{E.4}(effective) we know that the number of classes of binary forms of degree $n$ with fixed discriminant $v\prod{r_{i,j}^2}$ is finite and bounded by an integer depending only on $S$.\\
Since the cardinalities of the sets $\mathcal{R},\mathcal{V}$ are finite and depend only on $S$ and $K$, then the set of classes of $n$-tuples with good reduction outside $S$ has finite cardinality depending only on $S$ and $K$.\\
\end{proof}

\section{Cycles for rational maps}
Let $\Phi\colon\mathbb{P}_1\to\mathbb{P}_1$ a rational map defined over $K$ by $\Phi([x:y])=\left[F(x,y):G(x,y)\right]$, where $F,G\in R\left[x,y\right]$ have no common factor and are homogeneous of the same degree.
\begin{defin}\label{d1}
We say that a morphism $\Phi\colon\mathbb{P}_1\to\mathbb{P}_1$ defined over $K$ has good reduction at a prime ideal $\mathfrak{p}$ if there exists a morphism $\tilde{\Phi}\colon\mathbb{P}_1\to\mathbb{P}_1$ defined over $K(\mathfrak{p})$ with $\deg\Phi=\deg\tilde{\Phi}$ such that the following diagram
\begin{center}
$\begin{diagram}
\node{\ \mathbb{P}_{1,K}}\arrow{e,t}{\Phi} \arrow{s,l}{\widetilde{\phantom{\Phi}}}
\node{\ \mathbb{P}_{1,K}}\arrow{s,r}{\widetilde{\phantom{\Phi}}}\\
\node{\ \mathbb{P}_{1,K(\mathfrak{p})}}\arrow{e,t}{\widetilde{\Phi}}\node{\ \mathbb{P}_{1,K(\mathfrak{p})}}
\end{diagram}$\end{center} is commutative, where $\ \tilde{}\ $ denotes the reduction modulo $\mathfrak{p}$. Furthermore if $\Phi$ has good reduction at every prime ideal $\mathfrak{p}\notin S$, we say that it has good reduction outside $S$.
\end{defin}
We may assume that $F,G$ have coefficients in $R_{\mathfrak{p}}$ and that at least one coefficient is in $R_{\mathfrak{p}}^*$; therefore we obtain a  rational map, defined over $K(\mathfrak{p})$, as follows
\begin{displaymath} \tilde{\Phi}\colon\mathbb{P}_1\to\mathbb{P}_1;\ \ \tilde{\Phi}([x:y])=[\tilde{F}(x,y):\tilde{G}(x,y)],\end{displaymath} 
where $\tilde{F},\tilde{G}$ are the polynomials obtained by reduction modulo $\mathfrak{p}$ of the coefficients of $F$ and $G$.  This is the rational map which appears in the Definition \ref{d1}. Hence the rational map  $\Phi$, with its coefficients chosen as described above, has good reduction at the prime ideal $\mathfrak{p}$ if and only if  Res$(\tilde{F},\tilde{G})$ is non zero (since $F$ and $G$ have no common factors). 

If $H(t_1,\ldots,t_k)\in K[t_1,\ldots,t_k]$ is a non zero polynomial, following the notation of \cite{M.S.2}, we define $v_{\mathfrak{p}}(H)$ as
\begin{equation*} v_{\mathfrak{p}}(H)=v_{\mathfrak{p}}\left(\sum_{\substack{I}}a_{_I}t_1^{i_1}\cdots t_k^{i_k}\right)=\min_{\substack{I}}v_{\mathfrak{p}}(a_{_I})\end{equation*}
where the minimum is taken over all multi-indexes  $I=(i_i,\ldots,i_k)$. That is, $v_{\mathfrak{p}}(H)$ is the smallest valuation of the coefficients of $H$. For any family of polynomials $H_1,\ldots,H_m\in K[t_1,\ldots,t_k]$ we define $v_{\mathfrak{p}}(H_1,\ldots,H_m)$ to be the minimum of the $v_{\mathfrak{p}}(H_i)$, $i\in\{1,\ldots,m\}$.\\
Let $\Phi$ be defined as above and let Res$(F,G)$ be the resultant of homogeneous polynomials $F$ and $G$ of degree  $d$. We define Disc$(\Phi)$ to be the integral ideal of $R$ whose valuation at the prime ideal $\mathfrak{p}$ is given by
\begin{displaymath} v_{\mathfrak{p}}(\text{Disc}(\Phi))=v_{\mathfrak{p}}(\text{Res}(F,G))-2dv_{\mathfrak{p}}(F,G).\end{displaymath}
The definition is a good one by the properties of the resultant. Moreover choosing the coordinates of $F,G$ in $R_{\mathfrak{p}}$ and such that at least one coefficient is in $R_{\mathfrak{p}}^*$, we obtain that $v_{\mathfrak{p}}(\text{Disc}(\Phi))=v_{\mathfrak{p}}(\text{Res}(F,G))$. Therefore $\Phi$ has good reduction at the prime ideal $\mathfrak{p}$ if and only if $v_{\mathfrak{p}}(\text{Disc}(\Phi))=0$.
We will consider only cycles for rational maps which have good reduction outside $S$ and of degree $\geq 2$. These maps form a semigroup under composition on which the group ${\rm PGL}_2(R_S)$ acts by conjugation.
\begin{defin}\label{(S,n)} An ordered $n$-tuple of elements of $\mathbb{P}_1(K)$ which is a cycle for a rational map with good reduction outside $S$ will be called a $(S,n)$-cycle.\end{defin}

We will use the following elementary proposition included in the paper by Morton and Silverman \cite[Proposition 6.1]{M.S.2}.
\begin{prop}\label{p1}
Let $\Phi\colon\mathbb{P}_1\to\mathbb{P}_1$ be a rational map over $K$ which has good reduction at the prime ideal $\mathfrak{p}$ and let $P\in\mathbb{P}_1(K)$ be a periodic point for $\Phi$ with minimal period $n$. Then
\begin{description}

\item{(a)}\ \ \ \ \ $\delta_{\mathfrak{p}}\,(\Phi^i(P),\Phi^j(P))=\delta_{\mathfrak{p}}\,(\Phi^{i+k}(P),
\Phi^{j+k}(P))$ for every $i,j,k\in\mathbb{Z}$
\item{(b)} Let $i,j\in\mathbb{Z}$ such that $\gcd(i-j,n)=1$. Then 
\begin{center} $\delta_{\mathfrak{p}}\,(\Phi^i(P),\Phi^j(P))=\delta_{\mathfrak{p}}\,(\Phi(P),P)$.\end{center}
\end{description} \begin{flushright}$\square$\end{flushright}\end{prop}
This proposition states that, for every indices $i,k$, the ideals $\mathfrak{I}(P_0,P_i)$ and $\mathfrak{I}(P_k,P_{k+i})$ defined in (\ref{I_i}) are equal. Moreover if $k$ and $n$ are coprime, then the ideals $\mathfrak{I}(P_0,P_k)$ and $\mathfrak{I}(P_0,P_1)$ are equal.
\subsection{Proof of Theorem \ref{tp} in the case that $R_S$ is a P.I.D..}\label{s.s.pid}
Since we want to obtain a finiteness result about rational maps with good reduction outside $S$, without loss of generality in this section we can suppose that $R_S$ is a P.I.D.. Indeed each rational map which has good reduction outside $S$ has good reduction outside every set of places which contains $S$; hence we can enlarge $S$ so that $R_S$ becomes a P.I.D.. Note that the cardinality of a minimum enlarged set, with the above property, is bounded by $s+h_{R_S}-1$, where $h_{R_S}$ is the class number of $R_S$. In fact: if all prime ideal are principal, then $R_S$ is a P.I.D.; otherwise such a prime ideal is contained in an ideal class which is not the trivial one, so if we add this prime ideal to $S$ (obtaining a larger set $S^{\prime}$) we have that the new ring $R_{S^{\prime}}$ has class number $h_{R_{S^{\prime}}}< h_{R_S}$; by inductive method it results that to obtain a P.I.D. it suffices to add to $S$ a number of prime ideals $\ \leq h_{R_S}-1$.

In all part of this subsection we suppose that $R_S$ is a principal domain so we can adopt the convention that
any point $P_i\in\mathbb{P}_1(K)$ will be represented by coprime $S$-integral homogeneous coordinates $[x_i:y_i]$.
By this convention for all prime ideals $\mathfrak{p}\notin S$ and every points $P_1,P_2\in\mathbb{P}_1(K)$ it follows that $\delta_{\mathfrak{p}}(P_1,P_2)=v_{\mathfrak{p}}(x_1y_2-x_2y_1)$.
The next lemma states part \emph{(a)} of the last proposition in a form which will be useful in the sequel:

\begin{lemma}\label{o1} For every $(S,n)$-cycle $(P_0, P_1,\ldots,P_{n-1})$ and for every $i,j\in\mathbb{Z}$ there exist an $S$-unit $u_{j,j+i}\in R_S^\ast$ such that
\begin{equation}\label{f}(x_jy_{j+i}-x_{j+i}y_j)=(x_0y_i-x_iy_0)u_{j,j+i}.\end{equation}
\end{lemma}
\begin{proof} Proposition \ref{p1} asserts that, for every prime ideal $\mathfrak{p}\notin S$ and for every couple of indexes $i,j\in \mathbb{Z}$, we have that $\delta_{\mathfrak{p}}(P_{j},P_{j+i})=\delta_{\mathfrak{p}}(P_0,P_i)$, therefore the identity $v_{\mathfrak{p}}(x_jy_{j+i}-x_{j+i}y_j)=v_{\mathfrak{p}}(x_0y_j-x_jy_0)$ holds. So we have 
\begin{displaymath} u_{j,j+i}=\frac{x_jy_{j+i}-x_{j+i}y_j}{x_0y_i-x_iy_0}\in R_S^\ast.\end{displaymath}
\end{proof}

Another simple but important fact is the following:

\begin{lemma}\label{ld}
For every $(S,n)$-cycle $(P_0, P_1,\ldots,P_{n-1})$ and for every prime ideal $\mathfrak{p}\notin S$ the following properties hold: 
\begin{enumerate}
\item  for all indexes $j\in \{0,1,\ldots,n-1\}, i\nequiv0$ $({\rm mod}\ n)$ we have $\delta_{\mathfrak{p}}(P_j,P_{j+i})\geq\delta_{\mathfrak{p}}(P_0,P_1)$, or equivalently  
\begin{equation}\label{C_i} C_i\coloneqq \frac{x_0y_i-x_iy_0}{x_0y_1-x_1y_0}\in R_S \end{equation} and
\begin{equation}\label{C_i/G} x_jy_{j+i}-x_{j+i}y_j=C_iu_{j,j+i}(x_0y_1-x_1y_0),\end{equation}
where $u_{j,j+i}\in R_S^\ast$.

\item let $P_0=\left[x_0:y_0\right]$ and $P_1=\left[x_1:y_1\right]$ be the first and the second point of the $(S,n)$-cycle. The matrix $A\in {\rm GL}_2(K)$ \begin{equation}\label{A} A=\begin{pmatrix}
	\displaystyle\frac{-y_0}{x_0y_1-x_1y_0} & \displaystyle\frac{x_0}{x_0y_1-x_1y_0}\\\\
	\displaystyle\frac{y_1}{x_0y_1-x_1y_0} & \displaystyle\frac{-x_1}{x_0y_1-x_1y_0}
\end{pmatrix}\end{equation} maps the vector $(x_0,y_0)$ to $(0,1)$ and the vector $(x_1,y_1)$ to $(1,0)$. For any index $k$, if $(\bar{x}_k,\bar{y}_k)^t=A(x_k,y_k)^t$, then for every indexes $j\in \{0,1,\ldots,n-1\}, i\nequiv0$ $({\rm mod}\ n)$ the following identities are verified 
\begin{equation}\label{C_{j-i}} \bar{x}_j\bar{y}_{j+i}-\bar{x}_{j+i}\bar{y}_j=-C_{i}u_{j,j+i},\end{equation}

where $u_{j,j+i}$ is the $S$-unit defined in part 1. Furthermore for every index $k>1$
\begin{equation}\label{nc} (\bar{x}_k,\bar{y}_k)=\left(C_k,-C_{k-1}u_{1,k}\right).\end{equation}

\item if $i,j$ are coprime integers, then
\begin{equation}\label{min d} \min\{\delta_{\mathfrak{p}}(P_0,P_i),\delta_{\mathfrak{p}}(P_0,P_j)\}=\delta_{\mathfrak{p}}(P_0,P_1)\end{equation}
and  
\begin{equation}\label{mN}\min\{v_{\mathfrak{p}}(C_i),v_{\mathfrak{p}}(C_j)\}=0\end{equation}
for every prime ideal $\mathfrak{p}\notin S$.
\end{enumerate}\end{lemma}
\begin{proof} \emph{1.}
Since the $\mathfrak{p}$-adic distance satisfies the following triangle inequality \cite[Proposition 5.1]{M.S.2}: \begin{displaymath} \delta_{\mathfrak{p}}(P,T)\geq \min\{\delta_{\mathfrak{p}}(P,Q), \delta_{\mathfrak{p}}(Q,T)\}\ \ \text{ for every $P,Q,T\in\mathbb{P}_1(K)$},\end{displaymath}
it follows that\begin{displaymath} \delta_{\mathfrak{p}}(P_j,P_{j+i})\geq\min\{\delta_{\mathfrak{p}}(P_j,P_{j+1}),\ldots , \delta_{\mathfrak{p}}(P_{j+i-1},P_{j+i})\}=\delta_{\mathfrak{p}}(P_0,P_1).\end{displaymath}
for any index $i,j$. Thus, by the choice of coprime homogeneous coordinates for every points of $\mathbb{P}_1(K)$ and part (a) of Proposition \ref{p1} we have that $v_{\mathfrak{p}}(x_0y_i-x_iy_0)\geq v_{\mathfrak{p}}(x_0y_1-x_1y_0)$, so (\ref{C_i}) is proved. The identity (\ref{C_i/G}) follows from (\ref{f}) and (\ref{C_i}).

\emph{2}. Let $A$ be the matrix defined in (\ref{A}). We have that $A(x_0,y_0)^t=(0,1)^t$ and $A(x_1,y_1)^t=(1,0)^t$. Putting $(\bar{x}_k,\bar{y}_k)^t=A(x_k,y_k)^t$ for every index $k>1$, the Equation (\ref{C_i/G}) implies that for every indexes $j\in \{0,1,\ldots,n-1\}, i\nequiv0$ $({\rm mod}\ n)$
\begin{equation}\label{bar} \bar{x}_j\bar{y}_{j+i}-\bar{x}_{j+i}\bar{y}_j=-\frac{x_jy_{j+i}-x_{j+i}y_j}{x_0y_1-x_1y_0}=-C_{i}u_{j,j+i},\end{equation}
since \mbox{$\det(A)=-(x_0y_1-x_1y_0)^{-1}$}, which proves (\ref{C_{j-i}}). \\Considering  (\ref{bar}) with $j=0,i=k$ and $j=1, k=i+1$ we prove (\ref{nc}).

\emph{3}. There exist $c,d\in \mathbb{Z}$ such that $ci+dj=1$. By part $1$, the triangle inequality and Proposition \ref{p1} it is verified that
\begin{displaymath} \delta_{\mathfrak{p}}(P_0,P_1)\leq\min\{\delta_{\mathfrak{p}}(P_0,P_{c i}),\delta_{\mathfrak{p}}(P_0,P_{-dj})\}\leq\delta_{\mathfrak{p}}(P_{ci},P_{-dj})=\delta_{\mathfrak{p}}(P_0,P_1)\end{displaymath}
and 
\begin{displaymath} \delta_{\mathfrak{p}}(P_0,P_i)\leq\delta_{\mathfrak{p}}(P_0,P_{ci})\ ;\ \delta_{\mathfrak{p}}(P_0,P_j)\leq\delta_{\mathfrak{p}}(P_0,P_{dj})\end{displaymath}
so (\ref{min d}) follows .\\
Now suppose that $\delta_{\mathfrak{p}}(P_0,P_i)=\delta_{\mathfrak{p}}(P_0,P_1)$, 
\begin{align*}
&\delta_{\mathfrak{p}}(P_0,P_i)=v_{\mathfrak{p}}(x_0y_i-x_iy_0)
=v_{\mathfrak{p}}(C_i)+v_{\mathfrak{p}}(x_0y_1-x_1y_0)
=v_{\mathfrak{p}}(C_i)+\delta_{\mathfrak{p}}(P_0,P_1)\end{align*} therefore we have that $v_{\mathfrak{p}}(C_i)=0$ which proves (\ref{mN}).
\end{proof}

Lemma \ref{ld} states that, for any $(S,n)$-cycle $(P_0, P_1,\ldots,P_{n-1})$ and for every couple of indexes $j\in \{0,1,\ldots,n-1\}, i\nequiv0$ $({\rm mod}\ n)$, the ideal $\mathfrak{I}(P_0,P_1)$ divides the ideal $\mathfrak{I}(P_j,P_{j+i})$. The $S$-integer $C_i$ generates the ideal $\mathfrak{I}(P_0,P_i)\cdot\mathfrak{I}(P_0,P_1)^{-1}$. Moreover if $i,j$ are coprime, then the greatest common divisor of $\mathfrak{I}(P_0,P_i)$ and $\mathfrak{I}(P_0,P_j)$ is $\mathfrak{I}(P_0,P_1)$.

\begin{lemma}\label{dC} Let $(P_0, P_1,\ldots,P_{n-1})$ be a $(S,n)$-cycle and let $i,j\in\mathbb{Z}$; then 
\begin{equation}\label{L_{i,j}} L_{i,j}\coloneqq \frac{x_0y_{i\cdot j}-x_{i\cdot j}y_0}{x_0y_j-x_jy_0}\in R_S.\end{equation}
Moreover, let $i,j$ be fixed coprime integers. If for every $(S,n)$-cycle the set of principal ideals generated by the possible values of $L_{i,j}$ is finite, then also the set of principal ideals generated by possible values of $C_i$ is finite, where $C_i$ is the $S$-integer defined in (\ref{C_i}).\end{lemma}
\begin{proof}Let $(P_0, P_1,\ldots,P_{n-1})$ be a cycle for a rational map $\Phi$. In order to prove the first part we simply apply Lemma \ref{ld} to the cycle $(P_0,P_j,P_{2j},\ldots)$, relative to rational map $\Phi^j$.\\
Now suppose that $i,j\in \mathbb{Z}$ are coprime. By (\ref{L_{i,j}}) and (\ref{C_i}) 
\begin{displaymath} x_0y_{i\cdot j}-x_{i\cdot j}y_0=L_{i,j}C_j(x_0y_1-x_1y_0)\end{displaymath} and \begin{displaymath} x_0y_{i\cdot j}-x_{i\cdot j}y_0=L_{j,i}C_i(x_0y_1-x_1y_0),\end{displaymath}
therefore it follows that $L_{i,j}C_j=L_{j,i}C_i$.\\
From (\ref{mN}) in Lemma \ref{ld} we deduce that $v_{\mathfrak{p}}(C_i)\leq v_{\mathfrak{p}}(L_{i,j})$, for every $\mathfrak{p}\notin S$, so the finiteness of principal ideals generated by the possible values of $L_{i,j}$ gives the finiteness of principal ideals generated by the possible values of $C_i$.
\end{proof}

In the next proofs we will frequently use the fact that \mbox{$S$-unit} equations of the type
\begin{equation}\label{a_i} a_1x_1+a_2x_2+\ldots+a_nx_n=1,\end{equation}
where $a_i\in K^*$, have only a finite number of non-degenerate solutions \begin{displaymath}(x_1,x_2,\ldots,x_n)\in (R_S^\ast)^n.\end{displaymath} A solution is called non-degenerate if no subsum vanishes (i.e. $\sum_{i\in I}{a_ix_i}\neq0$ for every nonempty subset $I\subsetneq\{1,2,\ldots,n\}$).\\
In other words, the equation $X_1+X_2+\ldots+X_n=1$ has only a finite number of non-degenerate solutions $(X_1,\ldots,X_n)\in K^n$ with $v_{\mathfrak{p}}(X_i)$ fixed for every index $i$ and for every $\mathfrak{p}\notin S$. Actually we shall use this result only for $n=2$ and $n=3$.

This equation has been widely studied in the literature. For $n=2$, the finiteness of solutions of equation (\ref{a_i}) was proved by C.L. Siegel in a particular case. Later K. Mahler studied the case $K=\mathbb{Q}$ and generic finite set $S$. In 1960 S. Lang extended Mahler's result to arbitrary fields $K$ of characteristic 0 and solutions in any group $\Gamma\subset K^*$ of finite rank. A. Baker obtained effective results using his bound for linear forms in logarithms.  J.-H. Evertse \cite{E.1}, studying the case where $K$ is a number field of degree $d$ over $\mathbb{Q}$, found that the set of solutions has cardinality $\leq3\cdot7^{d+2s}$. This upper bound depends only on $s=\# S$ and $d=\left[K:\mathbb{Q}\right]$. Note that $S$ includes all archimedean places of $K$, hence $s\geq d/2$. In this way we obtain the upper bound $3\cdot7^{4s}$ depending only on $s$.

For general $n>2$, A. J. van der Poorten and H. P. Schlickewei in \cite{VDP.1} and J. H. Evertse in \cite{E.3} showed that the set of solutions is finite (non effective results). The best quantitative result is due to J. H. Evertse \cite{E.2}, who found the upper bound $2^{35n^4s}$, depending only on $s$ and  $n$, for the cardinality of the set of solutions.

Using in the proof of Theorem \ref{tp} the bounds found by Evertse, we could obtain a quantitative result, depending explixitely on the cardinality of $S$ and the class number of $R_S$.

\begin{lemma}\label{ellittiche}
Let $D,E\in K^*$ be fixed. Given the equation 
\begin{equation}\label{J_2}y^2=Du+Ev,\end{equation}
the set 
\begin{displaymath} \left\{[u:v:y^2]\in \mathbb{P}_2(K)\mid (u,v,y)\in R_S^\ast\times R_S^\ast\times R_S\ \text{is a solution of (\ref{J_2})}\right\}\end{displaymath}
is finite and depends only on $S$, $D$ and $E$. In particular the subset of $R_S^\ast$ defined by
\begin{equation}\label{U} \left\{\frac{u}{v}\mid \text{$u,v\in R_S^\ast$ satisfy (\ref{J_2}) for a suitable $y\in R_S$}\right\}\end{equation} is finite and depends only on $S$, $D$ and $E$. The same assertion is valid for the set of principal ideals of $R_S$ defined by
\begin{equation}\label{C}\{yR_S\mid y\ \text{satisfies (\ref{J_2}) for suitable}\ u,v\in R_S^\ast\}.\end{equation} Moreover the finiteness of the last set is valid also in the case $DE=0$.\end{lemma}
\begin{proof} Since $R_S^\ast$ has finite rank, there exists a finite set $W\subset R_S^\ast$, depending only on $S$, such that for every $u\in R_S^\ast$ there exist $\bar{u}\in R_S^\ast$ and $w\in W$ such that $u=w\bar{u}^6$. Let $y$ be an integer which satisfies (\ref{J_2}) for suitable $u,v\in R_S^\ast$; then there exist $\bar{u},\bar{v}\in R_S^\ast$ and $w_1,w_2\in W$ such that  $y^2=Dw_1\bar{u}^6+Ew_2\bar{v}^6$. Therefore, we deduce that $(\bar{u}^2/\bar{v}^2,y/\bar{v}^3)$ is an $S$-integral point on the elliptic curve defined by the equation $Y^2=Dw_1X^3+Ew_2$. By the finiteness of $S$-integral points of elliptic curve (Siegel's Theorem, see for example \cite[Chapter 7]{S.1}) and finiteness of the set $W$ we deduce that 
\begin{equation}\label{u/v} \frac{y^2}{v}=\frac{1}{w_2}\left(\frac{y}{\bar{v}^3}\right)^2\ {\rm and}\ \frac{u}{v}=\frac{w_1}{w_2}\left(\frac{\bar{u}^2}{\bar{v}^2}\right)^3\end{equation} 
can assume only a finite number of values depending only on $S$. If $DE=0$, e.g. $E=0$, then it is trivial that $y^2R_S=DR_S$. This concludes the proof.
\end{proof}

\begin{lemma}\label{coniche}
Let $D,E\in K$ be fixed. Given the equation 
\begin{equation}\label{j_2} y^2=D^2u^2+DuEv+E^2v^2,\end{equation}
the set of ideals of $R_S$ defined by \begin{displaymath}\{yR_S\mid (u,v,y)\in R_S^\ast\times R_S^\ast\times R_S\ \text{is a solution of (\ref{j_2})}\}\end{displaymath} is finite and depends only on $S$, $D$ and $E$.\end{lemma}
\begin{proof} If $DE=0$ the lemma is trivial.\\
Let $DE\neq0$. We can suppose that $D$ and $E$ are integers. If they are not, we can choose an integer $F$, depending only on $D$ and $E$, such that $FD,FE\in R_S$ and replace $y^2$ with $F^2y^2$ in (\ref{j_2}).\\
If $(u,v,y)\in R_S^\ast\times R_S^\ast\times R_S\ \text{is a solution of (\ref{j_2})}$ then 
\begin{displaymath} y^2=D^2u^2+DuEv+E^2v^2=(Du-\zeta Ev)(Du-\bar{\zeta}Ev),\end{displaymath}
where $\zeta$ and $\bar{\zeta}$ are the primitive third roots of unity.
The elements $(Du-\zeta Ev)$ and $(Du-\bar{\zeta}Ev)$ of the extension $K(\zeta)/K$ are integers with the property that \begin{displaymath}(Du-\zeta Ev)-(Du-\bar{\zeta}Ev)=(\bar{\zeta}-\zeta)Ev.\end{displaymath}
Let $T$ be the ring of algebraic integers of $K(\zeta)$ and let $\bar{S}$ be the finite set defined by all places of $K(\zeta)$ which lie over any place of $K$ included in $S$. We can enlarge $\bar{S}$ to a finite set of places such that $T_{\bar{S}}$ is a unique factorization domain and such that $(\bar{\zeta}-\zeta)E\in T_{\bar{S}}^*$. The $T_{\bar{S}}$-integers $(Du-\zeta Ev)$ and $(Du-\bar{\zeta} Ev)$ are coprime and therefore, after
multiplication by a unit, are squares in $T_{\bar{S}}$ and so they can be expressed in the form $w\bar{y}^2$ with $w\in T_{\bar{S}}^*$ and $\bar{y}\in T_{\bar{S}}$. Applying the last lemma to the equation $\bar{y}^2=Du/w-\zeta Ev/w$ we easily deduce that there exists a finite set $\mathcal{U}$, depending only on $T_{\bar{S}}$, such that for every $u,v\in T_{\bar{S}}^*$ which satisfy (\ref{j_2}) we have that $u/v\in \mathcal{U}$.\\
Since $R_S^\ast\subset T_{\bar{S}}^*$ the last statement is true also when we consider $u,v\in R_S^\ast$ and the conclusion follows.
\end{proof}

\begin{lemma}\label{B_02} Let $n\geq3$ and $C_2$ be the integer defined in (\ref{C_i}) of Lemma \ref{ld} and associated to a $(S,n)$-cycle. The set of principal ideals of $R_S$
\begin{displaymath}  \left\{C_2R_S\mid\ (P_0, P_1,\ldots,P_{n-1}) \ \text{is a $(S,n)$-cycle} \right\}\end{displaymath} is finite and depends only on $S$ and $K$.\end{lemma}
\begin{proof} We use the notation of Lemma \ref{ld}. Let $(P_0, P_1,\ldots,P_{n-1})$ be a $(S,n)$-cycle. 

Let us consider first the case that $n$ is an odd number. By Proposition \ref{p1}-(b) it follows that $\delta_{\mathfrak{p}}(P_0,P_1)$ and $\delta_{\mathfrak{p}}(P_0,P_2)$ are equal for every $(S,n)$-cycle and for every $\mathfrak{p}\notin S$; hence $v_{\mathfrak{p}}(x_0y_2-x_2y_0)=v_{\mathfrak{p}}(x_0y_1-x_1y_0)$, so from (\ref{C_i}) we deduce that $C_2=1$. This case is proved.

By Lemma \ref{ld}-Part 2, any $(S,n)$-cycle is mapped by the automorphism defined in (\ref{A}) to the following ordered $n$-tuple of vectors of $K^2$
\begin{displaymath} (0,1);(1,0);\ldots;(C_i,-C_{i-1}u_{1,i});\ldots;(C_{n-1},-C_{n-2}u_{1,n})\ .\end{displaymath}
Let $n\geq4$ be an even number. Suppose first that $3\nmid n$. By Proposition \ref{p1}-(b) we have that $\delta_{\mathfrak{p}}(P_0,P_3)=\delta_{\mathfrak{p}}(P_0,P_1)$, therefore we deduce that $C_3=1$. By the identity (\ref{C_{j-i}}) of Lemma \ref{ld} applied with $i=1$, $j=2$ and considering $C_1=1$ follows that $-C_2^2u_{1,3}+u_{1,2}=-u_{2,3}$; thus we obtain that $C_2$ satisfies
\begin{displaymath} C^2_2=\frac{u_{1,2}}{u_{1,3}}+\frac{u_{2,3}}{u_{1,3}}\ .\end{displaymath}
Lemma \ref{ellittiche}, applied with $u_{1,2}/u_{1,3}=u$, $u_{2,3}/u_{1,3}=v$ and $C_2=y$, proves this case.

Suppose now $n=2\cdot3^k\cdot m$ with $m>1$ and $3\nmid m$. For every $(P_0, P_1,\ldots,P_{n-1})$, the $n$-tuple $(P_0,P_{3^k},\ldots,P_{(2m-1)3^k})$ is a  $(S,2m)$-cycle and $2m\geq4$ is coprime with $3$. Let  
\begin{displaymath} L_{2,3^k}=\frac{x_0y_{2\cdot3^k}-x_{2\cdot3^k}y_0}{x_0y_{3^k}-x_{3^k}y_0}\ .\end{displaymath}
Applying to the above cycle the reasoning used in the previous case we obtain that, varying the possible $(S,n)$-cycle, the set of principal ideals of $R_S$ generated by $L_{2,3^k}$ is finite. Now we simply apply Lemma \ref{dC} with $i=2$ and $j=3^k$.

The last case is $n=2\cdot3^k$. We first reduce to the case $k=1$. If $k>1$ we consider the cycle $(P_0,P_{3^{k-1}},\ldots,P_{5\cdot3^{k-1}})$ which has length 6 and if the lemma holds in the case $n=6$, then one has the finiteness of ideals generated by $L_{2,3^k}$. Therefore, by Lemma \ref{dC} applied with $i=2$ and $j=3^k$, this case is proved.\\
Now we suppose that $n=6$. By Lemma \ref{ld}-Part 2, any cycle $(P_0,P_1,P_2,P_3,P_4,P_5)$, by the matrix defined in (\ref{A}), is sent to the following ordered $6$-tuple of vector of $K^2$
\begin{displaymath}(0,1);(1,0);(C_2,-u_{1,2});(C_3,-C_2u_{1,3});(C_4,-C_3u_{1,4});(C_5,-C_4u_{1,5})\end{displaymath}
with $u_{1,i}\in R_S^\ast$ for every $i\in\{2,3,4,5\}$.\\
By Lemma \ref{o1} the identities $C_4=C_2u_{0,4}$ and $C_5=u_{0,5}$ hold for suitable $u_{0,4},u_{0,5}\in R_S^\ast$. We rewrite the above $6$-tuple as
\begin{displaymath}(0,1);(1,0);(C_2,-u_{1,2});(C_3,-C_2u_{1,3});(C_2u_{0,4},-C_3u_{1,4});(u_{0,5},-C_2u_{0,4}u_{1,5})\end{displaymath}
Imposing the identity (\ref{C_{j-i}}) of Lemma \ref{ld} (after simple simplification) we obtain that $C_2,C_3$ satisfy the following system of equations:
\begin{equation}\label{s1} \left\{ \begin{array}{ll}
\vspace{1,3mm}C_3u_{1,4}-u_{1,2}u_{0,4}=u_{2,4}& \textrm{(a)}\\
\vspace{1,3mm}C_2^2u_{0,4}u_{1,5}-u_{1,2}u_{0,5}=C_3u_{2,5}& \textrm{(b)}\\
C_2^2u_{0,4}^2u_{1,5}-C_3u_{1,4}u_{0,5}=u_{4,5} & \textrm{(c)}
\end{array}\right.\end{equation}
The equation (a) is obtained from (\ref{C_{j-i}}) with $j=2$ and $i=2$, (b) with $j=2$ and $i=3$, (c) with $j=4$ and $i=1$.
From (a) it follows that
\begin{equation}\label{C_3} C_3=u_{1,2}u_{0,4}u_{1,4}^{-1}+u_{2,4}u_{1,4}^{-1}.\end{equation}
Now multiplying (b) by $-u_{0,4}$ and adding (c) we obtain
\begin{displaymath} u_{0,4}u_{1,2}u_{0,5}+C_3u_{2,5}u_{0,4}-C_3u_{1,4}u_{0,5}=u_{4,5},\end{displaymath}
replacing $C_3$ with the right term of (\ref{C_3}) in this last identity we obtain the following $S$-unit equation:
\begin{equation}\label{eq1} \frac{u_{1,2}u_{0,4}^2u_{2,5}}{u_{4,5}u_{1,4}}+\frac{u_{2,4}u_{0,4}u_{2,5}}{u_{4,5}u_{1,4}}-\frac{u_{2,4}u_{0,5}}{u_{4,5}}=1\end{equation}

Suppose that the equation (\ref{eq1}) does not admit vanishing subsums. By the $S$-units Theorem, there exist only finitely many possible values for the ratios
\begin{displaymath} \frac{u_{1,2}u_{0,4}^2u_{2,5}}{u_{4,5}u_{1,4}};\ \frac{u_{2,4}u_{0,4}u_{2,5}}{u_{4,5}u_{1,4}};\ \frac{u_{2,4}u_{0,5}}{u_{4,5}}.\end{displaymath}
From (\ref{C_3}) it follows that 
\begin{equation}\label{C_3.1} C_3=\frac{u_{4,5}}{u_{0,4}u_{2,5}}\left(\frac{u_{1,2}u_{0,4}^2u_{2,5}}{u_{4,5}u_{1,4}}+\frac{u_{2,4}u_{0,4}u_{2,5}}{u_{4,5}u_{1,4}}\right)\end{equation} therefore the set of principal ideals of $R_S$ generated by $C_3$ is finite. By the above equation (b) and a suitable application to Lemma  \ref{ellittiche}, also the set of ideals generate by $C_2$ is finite.

Suppose that in (\ref{eq1}) 
\begin{displaymath} \frac{u_{1,2}u_{0,4}^2u_{2,5}}{u_{4,5}u_{1,4}}+\frac{u_{2,4}u_{0,4}u_{2,5}}{u_{4,5}u_{1,4}}=0.\end{displaymath} 
In this case, by (\ref{C_3.1}), we have that $C_3=0$; this situation contradicts $n=6$.

Suppose that in (\ref{eq1}) 
\begin{displaymath} \frac{u_{1,2}u_{0,4}^2u_{2,5}}{u_{4,5}u_{1,4}}-\frac{u_{2,4}u_{0,5}}{u_{4,5}}=0;\end{displaymath}
which is equivalent to 
\begin{displaymath} \frac{u_{2,4}u_{0,4}u_{2,5}}{u_{4,5}u_{1,4}}=1.\end{displaymath}
From these last two identities it follows that 
\begin{equation}\label{sn2} u_{1,2}u_{0,4}u_{0,5}=\frac{u_{2,4}^2u_{0,5}^2}{u_{4,5}}.\end{equation}
Now multiplying (a) by $u_{0,5}$ and adding (c) we obtain
\begin{equation}\label{C_2^2} C_2^2=\frac{1}{u_{0,4}^2u_{1,5}}(u_{1,2}u_{0,4}u_{0,5}+u_{2,4}u_{0,5}+u_{4,5})\end{equation}
Replacing $u_{1,2}u_{0,4}u_{0,5}$ in (\ref{C_2^2}) with the right term of (\ref{sn2}) we obtain
\begin{displaymath} C_2^2=\frac{1}{u_{0,4}^2u_{1,5}u_{4,5}}(u_{2,4}^2u_{0,5}^2+u_{2,4}u_{0,5}u_{4,5}+u_{4,5}^2);\end{displaymath}
so applying in the suitable way Lemma \ref{coniche} we obtain the finiteness of the set of principal ideal of $R_S$ generated by $C_2$.
At last we consider the case 

\begin{equation}\label{sn3}  \frac{u_{2,4}u_{0,4}u_{2,5}}{u_{4,5}u_{1,4}}-\frac{u_{2,4}u_{0,5}}{u_{4,5}}=0\end{equation}
which is equivalent to 
\begin{displaymath} \frac{u_{1,2}u_{0,4}^2u_{2,5}}{u_{4,5}u_{1,4}}=1\end{displaymath} as well as 
\begin{equation}\label{sn3.1} \frac{u_{0,4}u_{2,5}}{u_{4,5}u_{1,4}}=\frac{1}{u_{1,2}u_{0,4}}.\end{equation}
Replacing in (\ref{sn3}) $(u_{0,4}u_{2,5})/(u_{4,5}u_{1,4})$ with the right term of (\ref{sn3.1}) we obtain
\begin{displaymath} \frac{u_{2,4}}{u_{1,2}u_{0,4}}-\frac{u_{2,4}u_{0,5}}{u_{4,5}}=0\end{displaymath}
which is equivalent to $u_{1,2}u_{0,4}u_{0,5}=u_{4,5}$.
From this last identity and (\ref{C_2^2}) we obtain
\begin{displaymath} C_2^2=\frac{1}{u_{0,4}^2u_{1,5}}(u_{2,4}u_{0,5}+2u_{4,5}).\end{displaymath}
By Lemma \ref{ellittiche} we have the finiteness of the set of principal ideal of $R_S$ generated by $C_2$.
This last case concludes the proof.
\end{proof}
\begin{coroll}\label{B_02^k} For every $l\in\mathbb{N}$, let $C_{2^l}$ be the integer defined in Lemma \ref{ld} associated to a $(S,n)$-tuple. Then the set of principal ideals of $R_S$  
\begin{displaymath} \left\{C_{2^l}R_S\mid\ (P_0, P_1,\ldots,P_{n-1}) \ \text{is a $(S,n)$-cycle}\right\}\end{displaymath} is finite and depends only on $l,S$ and $K$.\end{coroll}
\begin{proof} We use the same notation of Lemma \ref{ld}.\\
We prove the finiteness by induction on $l$. The case $l=1$ was already proved in Lemma \ref{B_02}. Suppose that the  statement is valid for $l-1$. Let
\begin{displaymath} L_{2,2^{l-1}}\coloneqq \frac{x_0y_{2^l}-x_{2^l}y_0}{x_0y_{2^{l-1}}-x_{2^{l-1}}y_0}\ .\end{displaymath} 
Applying Lemma \ref{B_02} to the cycle $(P_0,P_{2^{l-1}},P_{2^l},\ldots)$ we get the finiteness of the set of principal ideals generated by the possible values of $L_{2,2^{l-1}}$. \\
By Lemma \ref{ld} 
\begin{displaymath} C_{2^l}(x_0y_1-x_1y_0)=L_{2,2^{l-1}}C_{2^{l-1}}(x_0y_1-x_1y_0).\end{displaymath} Therefore, by inductive assumption on $C_{2^{l-1}}$, the conclusion follows.\end{proof}
\begin{lemma}\label{B_03} Let $C_3$ be the integer defined in Lemma \ref{ld} associated to a $(S,n)$-tuple. The set of principal ideals of $R_S$
\begin{displaymath}  \left\{C_3R_S\mid\ (P_0, P_1,\ldots,P_{n-1}) \ \text{is a $(S,n)$-cycle} \right\}\end{displaymath} is finite and depends only on $S$ and $K$.\end{lemma}
\begin{proof} We use the same notation of Lemma \ref{ld}. The statement is trivial for $n< 4$. \\For $n=4$ the lemma follows from Proposition \ref{p1} , since $n$ and $3$ are coprime. \\Let $n>4$. By Lemma \ref{ld}, any $(S,n)$-cycle is sent by the automorphism defined in (\ref{A}) to the following ordered $n$-tuple of vectors of $K^2$:
\begin{displaymath} (0,1);(1,0);(C_2,-u_{1,2});(C_3,-C_2u_{1,3});(C_4,-C_3u_{1,4});\ldots;(C_{n-1},-C_{n-2}u_{1,n}).\end{displaymath}
Thus by the identity in (\ref{C_{j-i}}), with $j=3,i=1$ it follows that
\begin{displaymath} -C_3^2u_{1,4}+C_4C_2u_{1,3}=u_{3,4}.\end{displaymath} 
Note that by Corollary \ref{B_02^k} we can choose a finite set $\mathscr{C}_2$ (the choice depends only on $S$) such that $C_2=D_2w_1$ and $C_4=D_4w_2$ with $D_2,D_4\in \mathscr{C}_2$ and $w_1,w_2$ suitable $S$-units. Then $C_3$ satisfies one of the finitely many equations 
\begin{equation}\label{eq1l9} C_3^2=D_4D_2\frac{w_1w_2u_{1,3}}{u_{1,4}}-\frac{u_{3,4}}{u_{1,4}}.\end{equation} 
Thus applying Lemma \ref{ellittiche} to equation (\ref{eq1l9}), with $w_1w_2u_{1,3}/u_{1,4}=u$, $u_{3,4}/u_{1,4}=v$ and $C_3=y$, the lemma is proved.
\end{proof}

\begin{oss}\label{st}\rm Note that the sets of ideals defined in Lemma \ref{ellittiche} and Lemma \ref{coniche} have cardinality bounded by a constant depending only on $S$, the coefficient $D,E$ and the field $K$ (see the proof Theorem D.8.3. in \cite{H.S.1} furthermore we can effectively determine these ideal sets).  So the cardinalities of the sets of ideals defined in Lemma \ref{B_02} and Lemma \ref{B_03} are bounded by a constant depending only on $S$ and $K$.\vspace{-6.8mm}\begin{flushright}$\square$\end{flushright}
\end{oss}
\medskip
\noindent Recall that the cross-ratio of four distinct points $P_1,P_2,P_3,P_4$ of $\mathbb{P}_1(K)$ is
\begin{displaymath} \varrho(P_1,P_2,P_3,P_4)=\frac{(x_1y_3-x_3y_1)(x_2y_4-x_4y_2)}
{(x_1y_2-x_2y_1)(x_3y_4-x_4y_3)}.
\end{displaymath}

\begin{proof}(\emph{Theorem 1}) In this proof we will use the notation of Lemma \ref{ld}. \\By Lemma \ref{ld}-part 1  it follows that, for every\ $(S,n)$-cycle\ $(P_0, P_1,\ldots,P_{n-1})$\ , the ideal $\mathfrak{I}_1=\prod_{\substack{\mathfrak{p}\notin S}}\mathfrak{p}^{\delta_{\mathfrak{p}}(P_0,P_1)}$ divides the ideal $\mathfrak{I}_i=\prod_{\substack{\mathfrak{p}\notin S}}\mathfrak{p}^{\delta_{\mathfrak{p}}(P_0,P_i)}$, for any index $i$. This proves that every fractional ideal 
\begin{displaymath} \mathfrak{I}_i\mathfrak{I}_1^{-1}=\prod_{\substack{\mathfrak{p}\notin S}}\mathfrak{p}^{\delta_{\mathfrak{p}}(P_0,P_i)-\delta_{\mathfrak{p}}(P_0,P_1)}\end{displaymath}
is actually an integral ideal of $R_S$.

The theorem is equivalent to proving the finiteness of the possible values for $\delta_{\mathfrak{p}}(P_0,P_i)-\delta_{\mathfrak{p}}(P_0,P_1)$ for any $(S,n)$-cycle and $\mathfrak{p}\notin S$ .

Let $(P_0, P_1,\ldots,P_{n-1})$ be a $(S,n)$-cycle. The matrix  $A$ defined in (\ref{A}) defines an element of ${\rm PGL}_2(K)$ which maps the ordered $n$-tuple $(P_0, P_1,\ldots,P_{n-1})$ to the $n$-tuple\ $(\bar{P_0},\bar{P_1},\ldots,\bar{P}_{n-1})$. We represent the points of the last $n$-tuple with coordinates as defined in (\ref{nc}) such that the equations (\ref{C_{j-i}}) are satisfied.\\
For every index $i$ the following identities hold:
\begin{equation}\label{d}\begin{split} \delta_{\mathfrak{p}}(\bar{P}_0,\bar{P}_i) & =v_{\mathfrak{p}}(\bar{x}_0\bar{y}_i-\bar{x}_i\bar{y}_0)\\
&=v_{\mathfrak{p}}(x_0y_i-x_iy_0)-v_{\mathfrak{p}}(x_0y_1-x_1y_0)\\
&=\delta_{\mathfrak{p}}(P_0,P_i)-\delta_{\mathfrak{p}}(P_0,P_1).\end{split}\end{equation}
Therefore by (\ref{d}) we get that for every $\mathfrak{p}\notin S$ the following statement is true:\\
\emph{the set 
\begin{displaymath} \Delta(\mathfrak{p},i)\coloneqq\left\{\delta_{\mathfrak{p}}(\bar{P}_0,\bar{P}_i)\mid (P_0, P_1,\ldots,P_{n-1}) \text{ is a $(S,n)$-cycle}\right\}\end{displaymath}
is finite depending only on $S$ and $K$ if and only if the set  \begin{displaymath} \left\{\delta_{\mathfrak{p}}(P_0,P_i)-\delta_{\mathfrak{p}}(P_0,P_1)\mid (P_0, P_1,\ldots,P_{n-1}) \text{ is a  $(S,n)$-cycle}\right\}\end{displaymath} is finite depending only on $S$ and $K$.}\\
Now we verify the finiteness of $\Delta(\mathfrak{p},i)$ for every index $i$ and for every $\mathfrak{p}\notin S$.\\
As remarked in the proof of Lemma \ref{B_03} we can choose a set $\mathscr{C}_2$ (the choice depends only on $S$ and $K$) such that for any possible value of $C_2$ we have that $C_2=D_2u$ with suitable $D_2\in \mathscr{C}_2$ and $u\in R_S^\ast$.\\
By (\ref{nc}) the point $\bar{P}_2$ is $[C_2:-u_{1,2}]$. Thus there exists a $S$--unit $u$ such that the matrix 
\begin{equation}\label{u} U=\begin{pmatrix}
	u & 0\\
	0 & -u_{1,2}
\end{pmatrix}\in{\rm PGL}_2(R_S) \end{equation} maps $\bar{P}_2$ to $[D_2:1]$ with $D_2\in\mathscr{C}_2$. Let $\tilde{P}_i=U(\bar{P}_i)=[\tilde{x}_i:\tilde{y}_i]$ for every index $i\geq3$; therefore, the automorphism $U$ maps the $n$-tuple $(\bar{P}_0,\bar{P}_1, \ldots,\bar{P}_{n-1})$ to
\begin{equation}\label{AU} ([0:1],[1:0],[D_2:1],\ldots,[\tilde{x}_i:\tilde{y}_i],\ldots,[\tilde{x}_{n-1}:\tilde{y}_{n-1}]),\end{equation}
where for every $\mathfrak{p}\notin S$ the $\mathfrak{p}$-adic distances are not changed.\\
It is clear that the theorem holds if $n<4$. Otherwise by the properties of the cross-ratio
\begin{displaymath} \varrho(\tilde{P}_0,\tilde{P}_1,\tilde{P}_2,\tilde{P}_3)+\varrho(\tilde{P}_0,\tilde{P}_1,\tilde{P}_3,\tilde{P}_2)=1,\end{displaymath}
we have that
\begin{equation}\label{birapporto} \frac{(\tilde{x}_0\tilde{y}_2-\tilde{x}_2\tilde{y}_0)(\tilde{x}_1\tilde{y}_3-\tilde{x}_3\tilde{y}_1)}{(\tilde{x}_0\tilde{y}_1-\tilde{x}_1\tilde{y}_0)(\tilde{x}_2\tilde{y}_3-\tilde{x}_3\tilde{y}_2)}-\frac{(\tilde{x}_0\tilde{y}_3-\tilde{x}_3\tilde{y}_0)(\tilde{x}_1\tilde{y}_2-\tilde{x}_2\tilde{y}_1)}{(\tilde{x}_0\tilde{y}_1-\tilde{x}_1\tilde{y}_0)(\tilde{x}_2\tilde{y}_3-\tilde{x}_3\tilde{y}_2)}=1.\end{equation}
By Lemma \ref{B_03} we can choose a set $\mathscr{C}_3$ (the choice depends only on $S$ and $K$), such that for any possible value of $C_3$ we have that $C_3=D_3w$ with suitable $D_3\in \mathscr{C}_3$ and $w\in R_S^\ast$.\\
The following identities hold:
\begin{align}\label{y_3} \tilde{y}_3&=\widetilde{D}_2v_{1,3}\\ \label{x_3}\tilde{x}_3&=D_3v_{0,3}\\ \tilde{x}_2\tilde{y}_3-\tilde{x}_3\tilde{y}_2&=v_{2,3},\end{align}
for suitable $v_{1,3},v_{0,3},v_{2,3}\in R_S^\ast$, $\widetilde{D}_2\in\mathscr{C}_2$ and $D_3\in\mathscr{C}_3$ .\\
Let us rewrite equation (\ref{birapporto}) as the $S$-unit equation in $v_{1,3}/v_{2,3},\ v_{0,3}/v_{2,3}$
\begin{equation}\label{birapporto1}-D_2\widetilde{D}_2\left(\frac{v_{1,3}}{v_{2,3}}\right)+D_3\left(\frac{v_{0,3}}{v_{2,3}}\right)=1\ .\end{equation}
Note that for any $D_2,\widetilde{D}_2\in\mathscr{C}_2,\ D_3\in\mathscr{C}_3$\ fixed, there are only finitely many solutions $(v_{1,3}/v_{2,3},v_{0,3}/v_{2,3})$. \\Hence, by (\ref{x_3}) and (\ref{y_3})
\begin{displaymath} [\tilde{x}_3:\tilde{y}_3]=\left[D_3\frac{v_{0,3}}{v_{2,3}}\frac{v_{2,3}}{v_{1,3}}:\widetilde{D}_2\right]\end{displaymath}
There are only finitely many equations of the type (\ref{birapporto1}) since the sets $\mathscr{C}_2$ and $\mathscr{C}_3$ are finite; therefore the finiteness of the possible values of $\tilde{P}_3$ follows, for every $(S,n)$-cycle $(P_0, P_1,\ldots,P_{n-1})$. By the same argument, increasing by $1$ each index in  (\ref{birapporto}) we obtain that the set of possible points $\tilde{P}_i$ is finite, thus it follows the finiteness of the values of  $\delta_{\mathfrak{p}}(\tilde{P}_0,\tilde{P}_i)$, for every index $i\in\{2,\ldots,n-1\}$ and for every $\mathfrak{p}\notin S$.\\
The proof now follows simply applying Corollary B in \cite{M.S.1}, which states that if $(P_0, P_1,\ldots,P_{n-1})$ is a $(S,n)$-cycle for rational map of degree $\geq2$, then \begin{displaymath} n\leq [12(s+2)\log(5s+10)]^{4[K:\mathbb{Q}]}\leq [12(s+2)\log(5s+10)]^{8s}.\end{displaymath} 
\end{proof}
By Remark \ref{st} the cardinality of the set $\mathbb{I}_S$ depends only on $S$ and $K$.

\subsection{Proofs in the general case.}
\begin{proof}[Proof of Theorem \ref{tp}]We deduce the general case of Theorem \ref{tp} from the particular case in which $R_S$ is a P.I.D., treated in subsection 3.1. For every $S$ such that $R_S$ is not a P.I.D. there exist infinitely many prime ideals of $R_S$ that are not principal (it is a direct consequence of unique factorization in prime ideals and the \emph{Chinese Remainder Theorem}). Therefore if $R_S$ is not a P.I.D., then there exist two disjoint finite sets $S_1$, $S_2$ of prime ideal such that $R_{S\cup S_1}$ and $R_{S\cup S_2}$ are P.I.D.. Let $(P_0, P_1,\ldots,P_{n-1})$ be a cycle for a rational map of degree $\geq2$ with good reduction outside $S$. For every index $i\in\{1,\ldots,n-1\}$, let $\mathfrak{I}_i$ be the ideals defined in (\ref{I_i}). For every prime ideal $\mathfrak{p}$ let $e_{\mathfrak{p}}$ be the exponent such that
\begin{displaymath}  \mathfrak{I}_i\mathfrak{I}_1^{-1}=\prod_{\substack{\mathfrak{p}\notin S}}\mathfrak{p}^{e_{\mathfrak{p}}}.\end{displaymath}
Applying Theorem \ref{tp} to $R_{S\cup S_1}$ we deduce that for every prime ideal $\mathfrak{p}\notin S\cup S_1$ (in particular for all $\mathfrak{p}\in S_2$) the value $e_{\mathfrak{p}}$ is non negative and bounded by a constant that depend only on $S$ and $K$. 
Now, applying Theorem \ref{tp} to $R_{S\cup S_2}$, the same statement holds for all prime ideal $\mathfrak{p}\in S_1$. In this way we have proved that there exist a constant $c$, which depends only on $S$ and $K$, such that for every prime ideal $\mathfrak{p}\notin S$, $0\leq e_{\mathfrak{p}}\leq c$. It is easy to see that for all prime ideal $\mathfrak{p}$, but finitely many (they are independent from the choice of $(S,n)$-cycle), $e_{\mathfrak{p}}$ is zero. This concludes the proof.
\end{proof}

\begin{proof}[Proof of Corollary \ref{P_0,P_1}] Fixing two consecutive points $P_0,P_1\in\mathbb{P}_1(K)$ of a $(S,n)$-cycle we set the ideal $\mathfrak{I}_1$ defined by (\ref{I_i}) and, by Theorem \ref{tp}, the set of possible ideals $\mathfrak{I}_i$ is finite and fixed . Hence the choice of two consecutive points of a $(S,n)$-cycle $(P_0, P_1,\ldots,P_{n-1})$ determines the finite set of possible values for $\delta_{\mathfrak{p}}(P_i,P_j)$, for any couple of points $P_i,P_j$. Applying the results of Birch and Merriman \cite{B.M.1}, with the same method used in the proof of Proposition \ref{fn}, we prove the corollary.
\end{proof}
\begin{proof}[Proof of Corollary \ref{N}] The proof of is contained in the proof of Theorem \ref{tp}. Indeed $\mathcal{N}$ is the set of the  $n$-tuples of type (\ref{AU}) which are obtained from any $(S,n)$-cycle $(P_0, P_1,\ldots,P_{n-1})$ under the action of the matrix $U\cdot A \in$ PGL$_2(K)$, where $A$ and $U$ are the matrices defined in (\ref{A}) and (\ref{u}), respectively. Since, for every $n$, the set of possible $n$-tuples of type (\ref{AU}) is finite (see the proof of Theorem \ref{tp}) and $n$ is bounded, then $\mathcal{N}$ is finite.
\end{proof}
\begin{proof}[Proof of Theorem \ref{tp2}]

Let $K=\mathbb{Q}$ and $S=\{|\cdot|_\infty;|\cdot|_2 \}$ so that $R_S=\mathbb{Z}[1/2]$. Let
\begin{displaymath} \mathcal{T}\coloneqq\left\{\left([u:u-1],[u-1:-1],[1:u]\right)\mid u\in R_S^\ast\right\}\end{displaymath} We will prove that the infinite set $\mathcal{T}$ is formed by $(S,3)$-cycles for suitable rational maps of degree equal to $4$ and that the set of primes defined by
\begin{displaymath}\{p\mid \delta_p(P_0,P_1)>0\text{ for some }(P_0,P_1,P_2)\in\mathcal{T}\}\end{displaymath}
is infinite. In particular the set of possible ideals $\mathfrak{I}_1$ is infinite. 
This will prove automatically the extension of Theorem \ref{tp2} to every number field $K$ and every choice of finite set $S$ containing all the  archimedean places of $K$ and the $2$-adic ones.

For every point $P=\left[x:y\right]\in\mathbb{P}_1(\mathbb{Q})$ we can choose coprime $S$-integral coordinates $(x,y)$. By this choice of coordinates, for every prime $p$ and for every couple of points $Q_1=\left[x_1:y_1\right]$, $Q_2=\left[x_2:y_2\right]$ of $\mathbb{P}_1(\mathbb{Q})$ the following identity holds
\begin{displaymath} \delta_p(Q_1,Q_2)=v_p(x_1y_2-x_2y_1).\end{displaymath}
Hence we have the following identity between ideals 
\begin{displaymath} (x_1y_2-x_2y_1)=\prod_{p\text{ prime}}p^{\delta_p(Q_1,Q_2)}.\end{displaymath}
To simplify the notation, to any rational map $\phi \colon\mathbb{P}_1\to\mathbb{P}_1$ defined over $\mathbb{Q}$ we associate, in the canonical way, the rational function  $\phi(z)\in \mathbb{Q}(z)$ by taking the pole of $z$ as the point at infinity $[1:0]$. In this way, a rational function $\phi(z)=N(z)/D(z)$ with $N,D\in\mathbb{Z}[z]$ coprime polynomials, has good reduction at a prime $p$ if $p$ does not divide the resultant of polynomials $F,D$ and $\tilde{\phi\phantom{.}}$, the rational function obtained from $\phi$ by reduction modulo $p$, has the same degree of $\phi$. Writing $\mathbb{P}_1(\mathbb{Q})=\mathbb{Q}\cup\{\infty\}$, we will shift from the homogeneous to the affine notation for points in $\mathbb{P}_1(\mathbb{Q})$ when necessary. So the point $[1:0]$ will correspond to $\infty$ and any other point $[x:y]$ will correspond to the rational number $x/y$.\\
Let $U(z)=(1-z)^{-1}$. Then for every $x/y\in \mathbb{Q}$ it follows that \begin{displaymath} U\left(\frac{x}{y}\right)=\frac{y}{y-x};\ \ U^2\left(\frac{x}{y}\right)=\frac{y-x}{-x}\ \ \ {\rm and}\ \ U^3\left(\frac{x}{y}\right)=\frac{x}{y}.\end{displaymath}
Indeed $U$ is the automorphism of $\mathbb{P}_1(\mathbb{Q})$ associated to the matrix $\begin{pmatrix} 0 & \hspace{-2.6mm}1\\-1 & 1 \end{pmatrix}$ of order $3$ in ${\rm PGL}_2(\mathbb{Z})$, so it has good reduction at any prime, since it is defined by a matrix in ${\rm SL}_2(\mathbb{Z})$. Note that every element of $\mathcal{T}$ is a cycle for $U$. Moreover $U$ admits the following cycles:
\begin{displaymath} 0\mapsto1\mapsto\infty\mapsto 0;\ \ \ -1\mapsto\frac{1}{2}\mapsto2\mapsto-1.\end{displaymath}
Define the degree three function $\Psi(z)\in\mathbb{Q}(z)$ by  
\begin{displaymath} \Psi(z)=\frac{(z+1)(2z-1)(z-2)}{2z(z-1)},\end{displaymath}
which has good reduction outside $S$ and satisfies $\Psi=\Psi\circ U$.\\
For every $S$-integer $u$, let us define $P_0=u/(u-1)$, $P_1=U(u/(u-1))=-(u-1)$ and $P_2=U^2(u/(u-1))=1/u$.
Since 
\begin{displaymath} \Psi(P_0)=\frac{-2u^3+3u^2+3u-2}{2u^2-2u}\end{displaymath}
we have that the automorphism defined by
\begin{displaymath} H(z)=\frac{(4u^2-4u)z-(-4u^3+6u^2+6u-4)}{2uz+4u^2+u-2}\end{displaymath} verifies $H(\Psi(P_0))=0$ and it belongs to ${\rm PGL}_2(\mathbb{Z}_S)$ if and only if $u\in R_S^\ast$. However the rational function
\begin{multline*}H\circ \Psi(z) =\frac{(2u^2-2u)z^3+(2u^3-6u^2+2)z^2+(-2u^3+6u-2)z+(2u^2-2u)}{uz^3+(2u^2-u-1)z^2+(-2u^2-2u+1)z+u}\end{multline*}
verifies 
\begin{equation}\label{0}H\circ \Psi(P_0)=H\circ \Psi(P_1)=H\circ \Psi(P_2)=0.\end{equation}
Define $\Psi_1(z)=z+H\circ \Psi(z)$.
\begin{lemma}\label{hg}Let $\phi\in\mathbb{Q}(z)$ be a rational map with good reduction outside $S$ such that $\phi(\infty)=\infty$, let $\begin{pmatrix} a & b\\u & c \end{pmatrix}\in {\rm GL}_2(R_S)$ and put
\begin{displaymath} T(z)=\frac{az+b}{uz+c}\in {\rm PGL}_2(R_S).\end{displaymath}
Then the rational function $z+T\circ \phi(z)$ has degree $\deg(\phi)+1$ and has good reduction outside $S$ if and only if $u\in R_S^\ast$.
\end{lemma}\begin{proof} Let $\phi(z)=N(z)/D(z)$ where $N,D\in \mathbb{Z}[z]$ are polynomials with no common factor. Since $\phi(\infty)=\infty$ we have that $\deg(\phi)=\deg(N)>\deg(D)$ and since $\phi$ has good reduction outside $S$ one has that the leading coefficient of $N$ is a $S$-unit and $N,D$ have no common factors modulo any prime $p\notin S$.
Also the rational function 
\begin{displaymath} T\circ \phi(z)=\frac{aN(z)+bD(z)}{uN(z)+cD(z)}\end{displaymath}
has good reduction outside $S$; therefore the polynomials $(aN(z)+bD(z)),(uN(z)+cD(z))$ have no common factor modulo any prime $p\notin S$ and have the same degree equal to $\deg(\phi)$; moreover the leading coefficient of $(uN(z)+cD(z))$ is a $S$-unit if and only if $u\in R_S^\ast$.
Now it is immediate to see that the rational function 
\begin{displaymath} z+T\circ \phi(z)=\frac{(uN(z)+cD(z))z+aN(z)+bD(z))}{(uN(z)+cD(z))}\end{displaymath}
has degree equal to $\deg(\phi)+1$ and has good reduction outside $S$ if and only if $u\in R_S^\ast$.\end{proof}

Since $\Psi(\infty)=\infty$, we can apply Lemma \ref{hg} with $\phi=\Psi$ and $T=H$ so $\Psi_1$ has good reduction outside $S$ if and only if $u\in R_S^\ast$. Moreover $\Psi_1$ has degree equal to $4$ and by (\ref{0}) follows that $\Psi_1(P_i)=P_i$ for every index $i\in\{0,1,2\}$. Thus we get that $(P_0,P_1,P_2)$ is a cycle for the rational map $\Phi=U\circ\Psi_1$ of degree 4 with good reduction outside $S$ (if $u\in R_S^\ast$).
The $3$-cycle \begin{displaymath}(P_0;P_1;P_2)=([u:u-1];[u-1:-1];[1:u])\in\mathcal{T},\end{displaymath}
with $u=2^n$ and $n\in\mathbb{N}$ gives $\mathfrak{I}_1=\mathfrak{I}_2=(2^{2n}-2^n+1)\cdot\mathbb{Z}[1/2]$ proving Theorem \ref{tp2}. Note also that the set of prime divisors of $(2^{2n}-2^n+1)$ for $n\in\mathbb{N}$ is infinite. This concludes the proof of Theorem \ref{tp2}.
\end{proof}
In the last proof we constructed a set of cycles of length 3. This is a consequence of the choice of the automorphism $U\in {\rm PGL}_2(\mathbb{Z})$ of order 3. However, for any $n$, with a suitable number field $K$ and a suitable finite set $S$ of places, it is possible to employ the same method as in the proof of Theorem \ref{tp2} starting with an automorphism of ${\rm PGL}_2(R_S)$ of order $n$. In this way it is possible to construct an infinite set of $(S,n)$-cycles which satisfies Theorem \ref{tp2}.

\bigskip

\noindent Jung Kyu CANCI\\Dipartimento di Matematica e Informatica\\Università degli Studi di Udine\\via delle Scienze, 206, 33100 Udine (ITALY).\\E-mail: {\sf canci@dimi.uniud.it}

\end{document}